\documentclass{birkjour}
\usepackage{bm}
\usepackage{mathrsfs}
\usepackage{amsfonts}
\usepackage{amssymb}
\usepackage{bm}
\usepackage{multirow,bigstrut}
\usepackage{threeparttable}
\usepackage{enumerate}
\usepackage{tikz}
\usepackage{graphicx}
\usepackage{stfloats}
\usepackage{array}
\usepackage{amsmath}

\usepackage{algorithm}
\usepackage{algorithmic}
\usepackage{setspace}

 \newtheorem{thm}{Theorem}[section]
 
 \newtheorem{lem}[thm]{Lemma}
 \newtheorem{prop}[thm]{Proposition}
 \theoremstyle{definition}
 \newtheorem{defn}[thm]{Definition}
 \theoremstyle{remark}
 \newtheorem{rem}[thm]{Remark}
 \newtheorem*{ex}{Example}
 \numberwithin{equation}{section}

 \newcommand{\Q}{\mathbb{H}}
\newcommand{\qi}{{\bf i}}
\newcommand{\qj}{{\bf j}}
\newcommand{\qk}{{\bf k}}
\newcommand{\bfu}{{\bm{\mu}}}

\newcommand{\e}{{\bf e}}
\newcommand{\sinc}{\mathrm{sinc}}
\newcommand{\Sc}{{\mathrm{Sc}}}
\newcommand{\Ve}{{\mathrm{Vec}}}

\newcommand{\bx} {\mathbf{x}}

\newcommand{\X} {\mathbf{X}}
\newcommand{\D} {\mathbf{D}}
\newcommand{\w} {\bm{\omega}}
\newcommand{\n} {\mathbf{n}}

\newcommand{\y} {\mathbf{y}}

\newcommand{\norm}[1]{\left\lVert#1\right\rVert}

\begin{document}

%-------------------------------------------------------------------------
% editorial commands: to be inserted by the editorial office
%
%\firstpage{1} \volume{228} \Copyrightyear{2004} \DOI{003-0001}
%\seriesextraline{This is the Concrete Title of this Book\br H.E. R and S.T.C. W, Eds.}
%\firstpage{1}

\title[Sampling Formulas Associated with QPSWFs]{Novel Sampling Formulas Associated with Quaternionic Prolate Spheroidal Wave functions}

%----------Author 1
\author[D. Cheng]{Dong Cheng}

\address{%
Department of Mathematics\\
University of Macau\\
Macao\\
 China}

\email{chengdong720@163.com}

%\thanks{This work was completed with the support of our\TeX-pert.}
%----------Author 2
\author[K. I. Kou]{Kit Ian Kou}
\thanks{Corresponding author: K. I. Kou}
\address{%
Department of Mathematics\\
University of Macau\\
Macao\\
 China}
\email{kikou@umac.mo}


\date{}

\begin{abstract}
The Whittaker-Shannon-Kotel'nikov (WSK) sampling theorem provides a reconstruction formula for the bandlimited signals. In this paper, a novel kind of the WSK sampling theorem is established by using the theory of quaternion reproducing kernel Hilbert spaces. This generalization is employed to obtain the novel sampling formulas for the bandlimited quaternion-valued signals. A special case of our result is to show that the 2D generalized prolate spheroidal wave signals obtained by Slepian can be used to achieve a sampling series of cube-bandlimited signals. The solutions of energy concentration problems in quaternion Fourier transform are also investigated.
\end{abstract}

%%% ----------------------------------------------------------------------
\maketitle
%%% ----------------------------------------------------------------------
%\tableofcontents
%%%%%%%%%%%%%%%%%%%%%%%%%%%%%%%%%%%%%%%%%%%%%%%%%%%%%%%%%%%%%%%%%%%%%%%%%%%%%%%%%%%%%%%%%%%%%%%%%%%%%%%
\section{Introduction}\label{S1}
 {T}he sampling theory, as one of the basic and fascinating topics in engineering sciences, is
crucial for reconstructing the continuous signals from the information collected at a series of discrete points without aliasing because it bridges the continuous physical signals and the discrete domain. After the celebrated Whittaker Shannon Kotel'nikov (WSK) sampling theorem established, there have been numerous proposals in the literature to generalize the classical WSK sampling expansions in various areas. The goal of this paper is to use the theory of quaternion reproducing kernel Hilbert method to obtain a generalization of WSK sampling for a general class of bandlimited quaternion-valued signals.

On the other hand, special functions \cite{bell1968special} such as Hermite and Laguerre functions have played an important role in classical analysis and mathematical physics. In a series of papers, Slepian {\it et al.} \cite{landau1961prolate1,landau1961prolate,landau1962prolate,slepian1964prolate} extensively investigated the remarkable properties of the prolate spheroidal wave functions (PSWFs) which are a class of special functions. For fixed $\tau$ and $\sigma$, the PSWFs of degree $n$ denoted by $\varphi_{n}$ constitute an
orthogonal basis of the space of $\sigma$-bandlimited signals with finite energy,
that is, for continuous finite energy signals whose Fourier transforms have support
in $[-\sigma,\sigma]$. They are also maximally concentrated on the interval $[-\tau,\tau]$ and depend
on parameters $\tau$ and $\sigma$.
PSWFs are characterized as the eigenfunctions of an integral operator with kernel arising from the sinc functions  $\sinc(x)=\frac{\sin \pi x}{\pi x}$:
\begin{equation}\label{1d-integral-eqn}
  \frac{\sigma}{\pi}\int^\tau_{-\tau}\varphi_{n}(x)\sinc \left(\frac{\sigma}{\pi}(y-x) \right)dx=
  \mu_{n}\varphi_{n}(y),~~~|y|\leq\tau.
\end{equation}
It has been shown that (\ref{1d-integral-eqn}) has  solutions in $L^2([-\tau,\tau])$ only for a discrete set of real positive values of $\mu_{n}$ say $\mu_{1}>\mu_{2}>...$ and that $\lim_{n\rightarrow \infty}\mu_{n}=0$. The corresponding
solutions, or eigenfunctions, $\varphi_{1,}(y),\varphi_{2}(y),...$ can be chosen to be real and orthogonal on $(-\tau,\tau)$.

 The variational problem that led to (\ref{1d-integral-eqn}) only requires that equation to hold for $|y|\leq\tau$. With $\varphi_{n}(x)$ on the left-hand side of
(\ref{1d-integral-eqn}) gives for $|x|\leq\tau$, however, the left is well defined for all $y$. We use this to extend the range of definition of the $\varphi_{n}$'s and so define
\begin{equation*}
  \varphi_{n}(y):=\frac{\sigma}{\pi\mu_{n}}\int^\tau_{-\tau}\varphi_{n}(x)
  \sinc \left(\frac{\sigma}{\pi}(y-x) \right)dx,~~~|y|\geq\tau.
\end{equation*}
The eigenfunctions $ \varphi_{n}$ are now defined for all $y$.  This leads to a dual orthogonality
\begin{align}
    &  \int^\tau_{-\tau}\varphi_{n}(x)\varphi_{m}(x)dx=\mu_{n}\delta_{mn},
    \label{dualorth1} \\
   &   \int^\infty_{-\infty}\varphi_{n}(x)\varphi_{m}(x)dx=\delta_{mn}. \label{dualorth2}
\end{align}

In \cite{zayed2007generalization}, Zayed proved that there are other systems of functions  possess similar properties to  those of prolate spheroidal wave functions. Moumni and Zayed \cite{moumni2014generalization} then extended the results to the higher dimension and derived a novel sampling formula for general class of bandlimited functions. This sampling formula \cite{moumni2014generalization} is a generalization of Walter and Shen's result \cite{walter2003sampling} on sampling with the PSWFs. In the present paper, we study the quaternionic prolate spheroidal  wave functions (QPSWFs), which refine and extend the PSWFs. The QPSWFs are ideally suited to study certain questions regarding the relationship between quaternion-valued signals and their Fourier transforms. We illustrate how to apply the QPSWFs for the quaternion  Fourier transform to analyze Slepian's energy concentration problem and sampling theory. We address all the above issues and explore some basic facts of the arising quaternion-valued function theory.

 Quaternion and quaternion Fourier transform  have already shown advantages over complex and classical Fourier transform within color image processing,  computer graphics, and robotics communities, for their modelling of rotation, orientation, and cross-information between multichannel data, see for instance, \cite{sangwine1996fourier,ell2007hypercomplex,ell2013quaternion}. We are motivated to develop a counterpart of the PSWFs in quaternion algebra. We apply the theory of reproducing-kernel Hilbert spaces and compact normal operators on Hilbert spaces, which are used in \cite{moumni2014generalization} to the quaternion algebra.  The contributions of this paper are summarized as follows.

 \begin{enumerate}
   \item    Due to the non-commutative property of the quaternion ring,   it is hard to study the spectral
property of compact normal operators in quaternion  Hilbert
space,   there is a big difference between the
complex case and the quaternionic case to consider this problem. In   Section \ref{S3}, we develop a spectral theorem for compact normal operators.  More importantly, the eigenvalues are sorted by norm in descending order.
   \item  We  introduce the quaternionic  prolate spheroidal   wave functions (QPSWFs) in Section \ref{S4}. Moreover, a series of elegant properties of QPSWFs are derived.
   \item Two sampling formulas (\ref{shannon sampling}) and (\ref{sampling with QPSWFs})  for bandlimited quaternion-valued signals  are obtained by using the QPSWFs.
   \item The maximum-energy-problem (Slepian's energy concentration problem) is also investigated.
 \end{enumerate}

The paper is organized as follows. In the next section, we collect some basic concepts in quaternion analysis. In Section \ref{S3}, we derive a spectral theorem for compact normal operators of quaternionic Hilbert spaces. In Section \ref{S4}, we introduce the QPSWFs and some of their properties. Moreover, these functions are used to obtain two sampling formulas of bandlimited quaternion-valued signals. Some examples are presented to demonstrate the results. More importantly, the maximum-energy-problem (Slepian's energy concentration problem) is also investigated.

\section{Preliminaries}\label{S2}
\subsection{Quaternion Algebra}
Throughout the paper, let
\begin{equation*}
  \Q :=\{q=q_0+\qi q_1+\qj q_2+\qk q_3|~q_0,q_1,q_2,q_3\in\mathbb{R}\},
\end{equation*}
be the {\em{Hamiltonian skew field of quaternions}}, where the elements $\qi$, $\qj$ and $\qk$ obey the Hamilton's multiplication rules:
\begin{equation*}
\qi^2=\qj^2=\qk^2=\qi\qj\qk=-1.
\end{equation*}
For every quaternion $q:=q_0+\underline{q}$, $\underline{q}:=\qi q_1+\qj q_2+\qk q_3$, the scalar and
vector parts of $q$, are  defined as $\Sc(q):=q_0$ and $\Ve(q):=\underline{q}$, respectively. If $q=\Ve(q)$, then $q$ is called pure imaginary quaternion.
The quaternion conjugate is defined by $\overline{q} :=q_0-\underline{q}=q_0-\qi q_1-\qj q_2-\qk q_3$, and the norm $|q|$ of $q$ is defined as
$|q|^2:={q\overline{q}}={\overline{q}q}=\sum_{m=0}^{3}{q_m^2}$.
Then we have
\begin{equation*}
  \overline{\overline{q}}=q,~~~\overline{p+q}=\overline{p}+\overline{q},~~~\overline{pq}=\overline{q}~\overline{p},~~~|pq|=|p||q|,~~~~\forall p,q\in\Q.
\end{equation*}
Using the conjugate and norm of $q$, one can define the inverse of $q\in\Q\backslash\{0\}$ as $q^{-1} :=\overline{q}/|q|^2$. The quaternion has subsets $\mathbb{C}_{\bm\mu}:=\{a+b {\bm\mu}:a,b\in\mathbb{R}, |\bm\mu|=1,  \bm \mu=\Ve(\bm{\mu}) \}$ . For each fixed unit pure imaginary quaternion $\bm\mu$, $\mathbb{C}_{\bm\mu}$ is isomorphic to the complex plane.
%%%%%%%%%%%%%%%%%%%%%%%%%%%%%%%%%%%%%%%%%%%%%%%%%%%%%%%%%%%%%%%%%%%%%

\subsection{Quaternion Module and Quaternionic Hilbert Space }

In order to state our results, we shall need some further notations.  The left quaternion module is similar with the right quaternion module, except that the quaternion ring acts on the left. The right quaternion module version of  all these facts can be found e.g. in (\cite{brackx1982clifford,ghiloni2013continuous}).

\begin{defn}[Quaternion module]\label{2.2.1}
Let $H$ be a left quaternion module, that is, $H$ consists of an abelian group with a left scalar multiplication $(q,u)\mapsto qu$ from $\Q\times H$ into $H$, such that for all $u,v\in H$ and $p,q\in\Q$
$$(p+q)u=pu+qu,~~p(u+v)=pu+pu,~~(pq)u=p(qu).$$
\end{defn}

\begin{defn}[Quaternionic pre-Hilbert space]\label{2.2.2}
A left quaternion module $H$ is called quaternionic pre-Hilbert space if there exists a quaternion-valued function (inner product) $(\cdot,\cdot):H\times H\rightarrow\Q$ with the following properties:
\begin{enumerate}
  \item $(u,v)=\overline{(v,u)}$;
  \item $(pu+qv,w)=p(u,w)+q(v,w)$;
  \item $(u,u)\in\mathbb{R}^+$ and $(u,u)=0$ if and only if $u=0$;
  \item $(0,v)=(u,0)=0$;
  \item $(u,pv+qw)=(u,v)\overline{p}+(u,w)\overline{q}$;
\end{enumerate}
where  $p,q\in\Q$ and $u,v,w\in H$.
\end{defn}

For each $u\in H$, putting $\|u\|^2={(u,u)}$, the Cauchy-Schwarz inequality and triangular inequality (see \cite{brackx1982clifford}) hold as $|(u,v)|^2 \leq (u,u)(v,v)$ and $\|u+v\| \leq \|u\| +\|v\|$.  The quaternionic pre-Hilbert space $H$ is said to be a quaternionic Hilbert space if it is complete under the norm $\|\cdot\|$. In what follows, by the notation $H$, we mean a (left) quaternionic Hilbert space.

For each $A\subset H$, define $A^\perp :=\{u\in H|~(u,v)=0,\forall v\in A\}$ and
$$\displaystyle \mathrm{Span}(A)=\left\{\sum_{k=1}^{n} q_ku_k:q_k\in\Q,u_k\in A,n\geq1\right \}. $$ Define $U(H):=\{u\in H,\|u\|=1\}$. $A$ is called an orthonormal set in $H$ if $A\subset U(H)$ and for any $u,v\in A$, $(u,v)=0$ for $u\neq v$.

\begin{thm}\label{2.2.1}
Let  $E$ be an orthonormal set in $H$. Then the following statements are equivalent.
\begin{enumerate}
  \item $E$ is a maximal orthonomal set ( i.e. if $E'$ is an orthonomal set such that $E\subset E'$, then $E'= E$ ).
  \item $E$ is total in $H$, that is, $\overline{\mathrm{Span}(E)}=H$.
  \item $E^\perp=\{0\}$.
  \item   $u=\sum_{z\in E} (u,z)z$ holds for every  $u\in H $;
  \item  $ (u,v)=\sum_{z\in E}(u,z)(z,v)$ holds for all $u,v\in H$.
  \item   $\|u\|=\sum_{z\in E} |(u,z)|^2$ holds for every $ u\in H$.
 \end{enumerate}
\end{thm}

\begin{thm} \label{preb}
Let  $A$ be a left $\Q$-linear subspace in $H$. Then the following assertions hold.
\begin{enumerate}
  \item $A^\perp$ is a  left $\Q$-linear closed subspace of $H$.
  \item If $A$ is closed, then $A=A^{\perp\perp}$ and $H=A\oplus A^\perp$, every $u\in H$ admits a unique decomposition $u=u_1+u_2$ with $u_1\in A$ and $u_2\in A^\perp$.
  \item If $A$ is  closed, calling $P_A(u)=u_1,u\in H$, we obtain that $P_A$ is a projection operator in $H$. Moreover $A=\mathrm{Range}(P_A)$, $A^\perp=\mathrm{Null}(P_A)$.
 \end{enumerate}
\end{thm}

A left $\Q$-linear operator is a map $\mathcal{T}:H\rightarrow H$ such that
\begin{equation*}
  \mathcal{T}(pu+qv)=p\mathcal{T}(u)+q\mathcal{T}(v)
\end{equation*}
if $p,q\in\Q$ and $u,v\in H$. Such an operator is called bounded if there exists a constant $c\geq0$ such that for all $u\in H$, $\|\mathcal{T}u\|\leq c\|u\|$.
As in the complex case, the norm of a bounded $\Q$-linear operator $T$ is defined by
\begin{equation*}
  \|\mathcal{T}\|:=\sup\left\{{\|\mathcal{T}u\|}:\|u\|\leq1\right\}=\sup\left\{{\|\mathcal{T}u\|}:\|u\|=1\right\}.
\end{equation*}
The set of bounded left $\Q$-linear operators is denoted by $\mathcal{B}(H)$.

\begin{prop}
Equip $\mathcal{B}(H)$  with the metric $\mathrm{Dist}(\mathcal{T}_1,\mathcal{T}_2)=\norm{\mathcal{T}_1-\mathcal{T}_2}$. Then $\mathcal{B}(H)$ is a complete metric space.
\end{prop}

 For every $\mathcal{T}\in\mathcal{B}(H)$, the Riesz representation theorem, as proposed in  \cite{brackx1982clifford},  guarantees that there exists a unique operator $\mathcal{T}^*\in\mathcal{B}(H)$, which is called the adjoint of $\mathcal{T}$, such that for all $u,v\in H$, $(\mathcal{T}u,v)=(u,\mathcal{T}^*v)$.
\begin{defn}\label{2.3}
Like the complex case, an operator $\mathcal{T}\in \mathcal{B}(H)$ is said to be
\begin{enumerate}
                  \item self-adjoint if $\mathcal{T}=\mathcal{T}^*$;
                  \item positive  if $\mathcal{T}$ is self-adjoint and $(\mathcal{T}u,u)\geq 0$ for every $u\in H$;
                   \item normal if  $\mathcal{T}\mathcal{T}^* =\mathcal{T}^*\mathcal{T}$;
                   \item compact if for every bounded set $B$ of $H$, $\overline{\mathcal{T}(B)}$ is a compact set of $H$.
               \end{enumerate}
\end{defn}
The set of all compact operators on $H$ is denoted by $\mathcal{B}_0(H)$. Clearly, $\mathcal{B}_0(H)\subset\mathcal{B}(H)$. Indeed, $\mathcal{B}_0(H)$ is a closed subset of $\mathcal{B}(H)$ (see Theorem 2 of  \cite{fashandi2013compact}) just like the complex case.

%%%%%%%%%%%%%%%%%%%%%%%%%%%%%%%%%%%%%%%%%%%%%%%%%%%%%%%%%%%%%%%%%%%%%%%%%%%%%%%%%%%%%%%%%%%%%%%%
\section{A Spectral Theorem for Compact Normal Operators in Quaternionic Hilbert Space}\label{S3}

Over the years, the spectral properties of operators in quaternionic Hilbert space have been studied (see e.g. \cite{fashandi2013compact,farenick2003spectral,ghiloni2014spectral,Alpay2016unbounded,Alpay2016unitary,
Colombo2016Schatten}).  In this part, we use  different approaches from  those used in aforementioned literature to derive a new spectral theorem for compact normal operators.

As in the complex case, define the eigenvalue $q$ of $\mathcal{T}\in \mathcal{B}(H)$ by
\begin{equation*}
  \mathcal{T}u=qu,~~~u\in H\setminus \{0\}.
\end{equation*}
However, as mentioned in  \cite{ghiloni2013continuous}, the eigenspace of $q$ cannot be a left $\Q$-linear subspace. If $\lambda\neq 0$, then $\lambda u$ is an eigenvector of $\lambda q\lambda^{-1}$ rather than $q$. The eigenvalues of compact normal operators are not necessary to be real. It is hard to consider the spectral  property of compact normal operators, because there is a big difference between the complex case and the quaternionic case to consider this problem.

Therefore, we  need to consider the entire similarity orbit $\theta(q)$ of $q$ (see \cite{farenick2003spectral}):
\begin{equation*}
  \theta(q):=\left\{\lambda q\lambda^{-1}:\lambda\in\Q\setminus \{0\}\right\}=\left\{\lambda q\overline{\lambda}:\lambda\in\Q,|\lambda|=1\right\}.
\end{equation*}
If $q\in\mathbb{R}$, then $\theta(q)$ contains only one element. In all other cases, $\theta(q)$ contains infinitely many elements. But the following lemma indicates that  only two of those are complex.
\begin{lem}[\cite{farenick2003spectral}]
If $q\in\Q$ is nonreal, then there is a nonreal $\lambda\in \mathbb{C}_{\bfu}$ such that $\theta(q)\cap\mathbb{C}_{\bfu}=\left\{\lambda,\overline{\lambda}\right\}$. In particular, if $\lambda\in \mathbb{C}_{\bfu}$, then $\theta(\lambda)\cap\mathbb{C}_{\bfu}=\left\{\lambda,\overline{\lambda}\right\}$.
\end{lem}

\begin{lem}\label{conjugacy-class}
If $\theta(q_1)\cap\theta(q_2)\neq\emptyset$, then $\theta(q_1)=\theta(q_2)$.
\end{lem}

For each $\theta(q)$, define the corresponding characteristic set as
\begin{equation*}
  H_q :=\left\{u\in H:\mathcal{T}u=pu,p\in \theta(q) \right\}.
\end{equation*}
Clearly, $H_q$ is a left $\Q$-linear subspace of $H$ when $q\in \mathbb{R}$. As $H_q$ is not always a left $\Q$-linear subspace of $H$, we call $H_q$ an eigenspace only when $H_q$ is   a left $\Q$-linear subspace of $H$.

\begin{rem}\label{leftORright}
A quaternion $\lambda$ is said to be a left (right) eigenvalue of quaternion matrix $A$ if there exists a vector $x\in \Q^{n\times 1}\setminus\{0\}$ such that
\begin{equation*}
   Ax=\lambda x,~ ~(Ax=x \lambda).
\end{equation*}

   In the present paper, the notion of left eigenvalue for left $\Q$-linear operator is different from the notion of left eigenvalue for quaternion matrix. The map $x\mapsto Ax$ is a right $\Q$-linear operator rather than a left $\Q$-linear operator.  In fact, the notion of left eigenvalue for left $\Q$-linear operator is actually similar with the notion of right eigenvalue for quaternion matrix. To see this, define a left operator $T$ as $y\mapsto yA$ for $y\in \Q^{1\times n}\setminus\{0\}$. Then a left eigenvalue of $T$ is the element $\lambda$ such that $Ty=\lambda y$ with $y\in\Q^{1\times n}\setminus\{0\}$.
It implies that $y A=\lambda y$. Taking the conjugate transpose on both sides of equality, we have
 \begin{equation*}
   A^\dagger y^\dagger =(yA)^\dagger=(\lambda y)^\dagger=y^\dagger \overline{\lambda},
 \end{equation*}
 where '$^\dagger$' is the conjugate transpose operation. Therefore $\overline{\lambda}$ is a right eigenvalue of $A^\dagger$.
\end{rem}

In the following some fundamental results of complex Hilbert space (see for instance  \cite{rudin1987real,bollabas1999linear}) are carried to quaternionic Hilbert space.

% Given the similarity between complex Hilbert space and quaternionic Hilbert space, some proofs of the following results are  omitted and the proofs of complex version can be found in \cite{rudin1987real,bollabas1999linear}. We only present the proofs of those which
%deserve an explicit proof.

\begin{lem}
$\mathcal{T}\in \mathcal{B}(H)$ is a compact operator if and only if for every bounded sequence $\{u_n\}\subset H$, the sequence $\{\mathcal{T}u_n\}$ has a convergent subsequence.
\end{lem}

\begin{prop}
Let $\mathcal{T}\in \mathcal{B}(H)$  be a normal operator. Then the following assertions hold.
\begin{enumerate}
\item $\|\mathcal{T}u\|=\|\mathcal{T}^*u\|$ for every $u\in H$.
\item $\mathcal{T}u=qu$  if and only if $\mathcal{T}^*u=\overline{q}u$, where $u\in H$ and $q\in\Q$.
\item If $q_1$ and $q_2$ are eigenvalues of $\mathcal{T}$ and $\theta(q_1)\neq\theta(q_2)$, then $H_{q_1}\perp H_{q_2}$.
\end{enumerate}
\end{prop}

\begin{proof}
Statement {1}: Since $\mathcal{T}^*\mathcal{T}=\mathcal{T}\mathcal{T}^*$, then $$\|\mathcal{T}u\|^2=(\mathcal{T}u,\mathcal{T}u)=(\mathcal{T}^*\mathcal{T}u,u)
=(\mathcal{T}\mathcal{T}^*u,u)=(\mathcal{T}^*u,\mathcal{T}^*u)=\|\mathcal{T}^*u\|^2.$$

Statement {2}: It can be easily seen from $(\mathcal{T}u-qu,\mathcal{T}u-qu)=(\mathcal{T}^*u-\overline{q}u,\mathcal{T}^*u-\overline{q}u)$.

Statement {3}: If $\mathcal{T}u=\lambda_1 u~(\lambda_1\in \theta(q_1))$ and $\mathcal{T}v=\lambda_2 v~(\lambda_2\in \theta(q_2))$ , we have $\mathcal{T}^*u=\overline{\lambda_1} u$ and  $\mathcal{T}^*v=\overline{\lambda_2} v$ by statement {2}. Therefore
\begin{equation*}
  \lambda_1(u,v)=(\mathcal{T}u,v)=(u,\mathcal{T}^*v)=(u,\overline{\lambda_2}v)=(u,v)\lambda_2.
\end{equation*}
By Lemma \ref{conjugacy-class}, $\theta(q_1)\cap\theta(q_2)=\emptyset$, thus $(u,v)=0$.
\end{proof}

\begin{prop}\label{prop-compact-self-adjoint-operator}
Let $\mathcal{T}\in \mathcal{B}(H)$  be a compact self-adjoint operator. Then the following assertions hold.
\begin{enumerate}
  \item  $\|\mathcal{T}\|=\sup \left\{ |(\mathcal{T}u,u)|: {\|u\|\leq 1}\right\}$.
  \item  Every left eigenvalue of $\mathcal{T}$ is real-valued.
  \item  $\mathcal{T}$ has an eigenvalue of absolute value $\|\mathcal{T}\|$. Moreover, if $\mathcal{T}$ is  positive, then $\norm{\mathcal{T}}$ is an eigenvalue of $\mathcal{T}$.
\end{enumerate}
\end{prop}
\begin{proof}
Statement {1}: Let $M=\sup \left\{ |(\mathcal{T}u,u)|: {\|u\|\leq 1}\right\}$. For each $\|u\|\leq 1$, $$|(\mathcal{T}u,u)|\leq \|\mathcal{T}u\| \|u\|\leq \|\mathcal{T} \| \|u\|^2 \leq\|\mathcal{T}\|.$$
Thus $M\leq \|\mathcal{T}\|$. For each $\|u\|=1$, let $v$ be the element such that $\mathcal{T}u=\|\mathcal{T}u\|v$. Note that $|(\mathcal{T}x,x)|\leq M \|x\|^2$ for every $x\in H$ and
\begin{equation*}
  \|\mathcal{T}u\|=(\|\mathcal{T}u\|,v)=(\mathcal{T}u,v)=\overline{(v,\mathcal{T}u)}=(v,\mathcal{T}u)
  =(\mathcal{T}v,u)
\end{equation*}
for $\mathcal{T}=\mathcal{T}^*$. Therefore
\begin{equation*}
   \begin{split}
      4 \|\mathcal{T}u\| & =2(\mathcal{T}u,v)+2 (\mathcal{T}v,u)\\
       & =(\mathcal{T}u+\mathcal{T}v,u+v)-(\mathcal{T}u-\mathcal{T}v,u-v)\\
       &\leq M(\|u+v\|^2+\|u-v\|^2)\leq M(\|u\|^2+\|v\|^2)=4M.
   \end{split}
\end{equation*}
It follows that $\|\mathcal{T}\|=\sup\{\|\mathcal{T}u\|:  {\|u\|= 1}\}\leq M$.

Statement {2}: Suppose that $\lambda \in \Q$ is a left eigenvalue of $\mathcal{T}$. Then there is a $v\in H$ with $\|v\|=1$ satisfying $\mathcal{T}v=\lambda v$. Since $\mathcal{T}=\mathcal{T}^*$, then
\begin{equation*}
  \lambda=(\lambda v,v)=(\mathcal{T}v,v)=(v,\mathcal{T}v)=(v,\lambda v)=(v,v)\overline{\lambda}=\overline{\lambda}
\end{equation*}
is real-valued.  Also, $(\mathcal{T}u,u)=(u,\mathcal{T}u)=\overline{(\mathcal{T}u,u)}$ implies that $(\mathcal{T}u,u)$ is real-valued for every $u\in H$.

Statement {3}: Set $M_1=\sup\{(\mathcal{T}u,u): \|u\|=1\}$ and $M_2=\inf\{(\mathcal{T}u,u): \|u\|=1\}$. Without loss of generality, we assume that $M_1\geq M_2$. Then by Statement {2}, we have $\|\mathcal{T}\|=M_1$. From definition of $M_1$, there exists a sequence $\{u_n\}$ with $\|u_n\|=1$ satisfying $\lim_{n\to \infty} (\mathcal{T}u_n,u_n)=M_1$. Since $\mathcal{T}$ is compact, there is a subsequence $\{u_{n_k}\}$ such that  $\{\mathcal{T}u_{n_k}\}$ converging to an element of $H$, say $\lim_{k\to \infty} \mathcal{T}u_{n_k}=v_0 $. Note that
\begin{equation*}
  |(\mathcal{T}u_{n_k},u_{n_k})|\leq \|\mathcal{T}u_{n_k}\|\|u_{n_k}\| =\|\mathcal{T}u_{n_k}\|\leq\|T\|.
\end{equation*}
It follows that $\|\mathcal{T}\|\geq \|v_0\|\geq|M_1|$ and therefore $\|\mathcal{T}\|= \|v_0\|=|M_1|$.
Note that $M_1$ is real-valued and
\begin{equation*}
  \|\mathcal{T}u_{n_k}-M_1  u_{n_k}\|=\|\mathcal{T}u_{n_k}\|^2-2 M_1(\mathcal{T}u_{n_k},u_{n_k})+M_1^2\|u_{n_k}\|^2.
\end{equation*}
Thus $\lim_{k\to \infty} M_1  u_{n_k} =\lim_{k\to \infty} \mathcal{T}u_{n_k}=v_0$. Set $x_0=M_1^{-1}v_0$, then
\begin{equation*}
   \mathcal{T}x_0=\lim_{k\to \infty}\mathcal{T}u_{n_k}=v_0=M_1 x_0
\end{equation*}
Therefore $M_1$ is an eigenvalue of $\mathcal{T}$.
\end{proof}

Since every eigenvalue of compact self-adjoint operator is real. So it is
easy to obtain the following theorem.
\begin{thm}\label{spectal-theorem-self-joint-operator}
 Let $H$ be a left quaternionic Hilbert space and $\mathcal{T}\in \mathcal{B}(H)$ be a compact self-joint operator. Then there is a (possibly finite) sequence $\{\lambda_k\}$ of real numbers satisfying $|\lambda_1|\geq|\lambda_2|\geq\cdot\cdot\cdot$ and a sequence $\{\xi_k\}\subset U(H)$
such that
\begin{enumerate}
  \item $\{\lambda_k\}$ are eigenvalues of $\mathcal{T}$ and $\lambda_k\rightarrow 0$ if  $\{\lambda_k\}$ is infinite;
  %\item if $H_k=\left\{u\in H:Tu=pu,p\in \theta(\lambda_k) \right\}$, then  $\dim H_k<\infty$;
  %\item if $H_k=\left\{u\in H:Tu=pu,p\in \theta(\lambda_k) \right\}$,then $H_k\perp H_l$ for $k\neq l$;
  \item $\mathcal{T}\xi_k=\lambda_k\xi_k$ for every $k$. For every $u\in H$, we have
  \begin{equation*}
    \mathcal{T}u=\sum_{k=1}^\infty (u,\xi_k)\lambda_k\xi_k;
  \end{equation*}
  %\item if $P_k$ is the orthogonal projection from $H$ into $H_k$, then
  %\begin{equation*}
   % T=\sum_{k=1}^\infty \lambda_kP_k;
  %\end{equation*}
  \item  if $0$ is not an eigenvalue of $\mathcal{T}$, then $\{\xi_k\}$ is a orthonormal basis of $H$;
  \item   $\dim H_{\lambda_k}<\infty$ for $\lambda_k\neq0$.
\end{enumerate}
\end{thm}

\begin{proof}
If $\mathcal{T}=0$, we are done. Suppose that $\mathcal{T}\neq0$, then $\|\mathcal{T}\|\neq0$. By Proposition \ref{prop-compact-self-adjoint-operator}, there exists a pair $(\lambda_1,\xi_1)\in \mathbb{R}\times U(H)$  satisfying $|\lambda_1|=\|\mathcal{T}\|$ such that $\mathcal{T}\xi_1=\lambda_1\xi_1$. Let $W_1=\mathrm{Span}\{\xi_1\}$, $V_1=\{u\in H:u\perp \xi_1\}$, then $W_1$, $V_1$ are left $\Q$-linear closed  subspace of $H$ and $H=W_1\oplus V_1$. Moreover, $\mathcal{T}(V_1)\subset V_1$, since if $v\in V_1$, we have
$$(\mathcal{T}v,\xi_1)=(v,\mathcal{T}^*\xi_1)=(v,\overline{\lambda_1}\xi_1)=(v,\xi_1)\lambda_1=0.$$

$V_1$ is still a left quaternionic Hilbert space. Define $\mathcal{T}_1=\mathcal{T}|_{V_1}$, then $\mathcal{T}_1$ a compact self-joint operator on $V_1$. If $\mathcal{T}_1=0$, then for any $u\in H$, we have $u=w_1+v_1$, where $w_1\in W_1$, $v_1\in V_1$. Thus
\begin{equation*}
  \mathcal{T}u = \mathcal{T}w_1+\mathcal{T}v_1= \mathcal{T}(w_1,\xi_1)\xi_1+0=\mathcal{T}(u,\xi_1)\xi_1
   = (u,\xi_1)\mathcal{T}\xi_1=(u,\xi_1)\lambda_1\xi_1.
\end{equation*}

If $\mathcal{T}_1\neq0$, then  there exists a pair $(\lambda_2,\xi_2)\in \mathbb{R}\times U(V_1)$  satisfying $|\lambda_2|=\|\mathcal{T}_1\|$ such that $\mathcal{T}\xi_2=\lambda_2\xi_2$. Let $W_2=\mathrm{Span}\{\xi_1, \xi_2\}$, $V_2=\{u\in H:u\perp W_2\}$, then $W_2$, $V_2$ are left $\Q$-linear closed  subspace of $H$ and $H=W_2\oplus V_2$. Furthermore,  $\mathcal{T}(V_2)\subset V_2$. Let $\mathcal{T}_2=\mathcal{T}_1|_{V_2}=0$, then for any $u\in H$, we have
$$\mathcal{T}u=(u,\xi_1)\lambda_1\xi_1+(u,\xi_2)\lambda_2\xi_2.$$

 If $\mathcal{T}_2\neq0$, continue the above procedure. If there is a $n\in\mathbb{N}$ such that $\mathcal{T}_n=0$, then
  \begin{equation*}
    \mathcal{T}u=\sum_{k=1}^n  (u,\xi_k)\lambda_k\xi_k.
  \end{equation*}
  Otherwise, ${\lambda_k}$ is  infinite, and $|\lambda_1|\geq|\lambda_2|\geq\cdot\cdot\cdot$. From the definitions of $W_n$ and $V_n$, we have $\xi_k \perp \xi_l$ for $k\neq l$. Furthermore, $\lambda_k\rightarrow 0$ as $k\rightarrow \infty$. If not, there is a subsequence $\{\lambda_{s_m}\}$ such that $|\lambda_{s_m}|\geq \delta>0$. Then
  \begin{equation*}
        \|\mathcal{T}\xi_{s_m}- \mathcal{T}\xi_{s_n}\|^2 = \|\mathcal{T}\xi_{s_m}\|^2+\|\mathcal{T}\xi_{s_n}\|^2
       = |\lambda_{s_m}|^2+|\lambda_{s_n}|^2\geq 2\delta^2>0,
  \end{equation*}
   which is contradict with the compactness of $\mathcal{T}$. Let $W_\infty=\overline{\mathrm{Span}\{\xi_k:k\geq1\}}$ and $V_\infty=\{u\in H:u\perp W_\infty\}\subset V_k~(k=1,2,\cdot\cdot\cdot)$. If $u\in V_\infty$, then $$|(\mathcal{T}u,u)|=|({\mathcal{T}_ku,u})|\leq \|\mathcal{T}_k\|\|u\|^2=|\lambda_k|\|u\|^2\rightarrow 0,$$
   as $k\rightarrow\infty$. Therefore, $(\mathcal{T}u,u)=0$ for every $u\in V_\infty$. Given $u\in U(V_\infty)$, let $v\in U(V_\infty)$ be such that $\mathcal{T} u=\|\mathcal{T}u\|v$. Then $(\mathcal{T}u,v)=(u,\mathcal{T}v)=\|\mathcal{T}u\|$ and thus
\begin{equation*}
 \|\mathcal{T}u\|=(\mathcal{T}u,v) =
\frac{1}{4}\left[(\mathcal{T}(u+v),u+v)-(\mathcal{T}(u-v),u-v)\right]=0.
\end{equation*}
   Hence $\mathcal{T}|_{V_\infty}=0$. Therefore, for every $u\in H$, let $u=w_\infty+v_\infty$, where $w_\infty\in W_\infty$ and $v_\infty\in V_\infty$, we have
   \begin{equation*}
     \begin{split}
     \mathcal{T}u =& \mathcal{T}w_\infty+\mathcal{T }v_\infty
       = \mathcal{T}(\sum_{k=1}^\infty (w_\infty,\xi_k)\xi_k)+0 \\
       =&  \sum_{k=1}^\infty (w_\infty,\xi_k)\mathcal{T}\xi_k
       = \sum_{k=1}^\infty (w_\infty,\xi_k)\lambda_k\xi_k \\
       =&   \sum_{k=1}^\infty (u,\xi_k)\lambda_k\xi_k.
     \end{split}
   \end{equation*}
    The third  equality holds for the boundness of $\mathcal{T}$. If $0$ is not an eigenvalue of $\mathcal{T}$, then $H=W_\infty$ and $\{\xi_k\}$ is a orthonormal basis of $H$. $\dim H_{\lambda_k}<\infty$ for $\lambda_k\neq0$ as $\mathcal{T}$ is compact.
\end{proof}

\begin{thm}\label{spectal-theorem-positive-operator}
Let $H$ be a left quaternionic Hilbert space and $\mathcal{T}\in \mathcal{B}(H)$ be a compact positive (self-joint) operator. Then there is a (possibly finite) sequence $\{\mu_k\}$ of nonnegative real numbers satisfying $\mu_1>\mu_2>\cdot\cdot\cdot$ %and a sequence $\{H_k\}$ of (left) $\Q$-linear subspaces of $H$
such that
\begin{enumerate}
  \item  $\{\mu_k\}$ are eigenvalues of $\mathcal{T}$ and $\mu_k\rightarrow 0$ if  $\{\mu_k\}$ is infinite;
  %\item if $H_k=\left\{u\in H:Tu=pu,p\in \theta(\lambda_k) \right\}$, then  $\dim H_k<\infty$;
  %\item if $H_k=\left\{u\in H:Tu=pu,p\in \theta(\lambda_k) \right\}$,then $H_k\perp H_l$ for $k\neq l$;
  %\item if $T\xi_k=\lambda_k\xi_k$ and $\|\xi_k\|=1$ for every $k$. Then for every $u\in H$, we have
  %\begin{equation*}
   % Tu=\sum_{k=1}^\infty (u,\xi_k)\lambda_k\xi_k;
  %\end{equation*}
  %\item if $P_k$ is the orthogonal projection from $H$ into $H_k$, then
  %\begin{equation*}
   % T=\sum_{k=1}^\infty \lambda_kP_k;
  %\end{equation*}
  %\item if $0$ is not a eigenvalue of $T$, then $\{\xi_k\}$ is a orthonormal basis of $H$;
  \item  if $H_k=\left\{u\in H:\mathcal{T}u=pu,p\in \theta(\mu_k) \right\}$,then $H_k\perp H_l$ for $k\neq l$ and $\dim H_k<\infty$ for $\mu_k\neq0$.
\end{enumerate}
\end{thm}

An exposition of the spectral theory of normal matrices with quaternion entries was presented in \cite{farenick2003spectral} by Farenick and Pidkowich. They also obtained a spectral theorem for compact operators in $n$-dimensional quaternionic Hilbert spaces by establishing relations between compact operators and quaternion normal matrices.

\begin{thm} \label{spectral-theorem-for-n-dimension}
Assume that $H$ is an n-dimensional quaternionic Hilbert space.Then an operator $\mathcal{T}:H\rightarrow H$ is normal if and only if there exists an orthonormal set $E=\{\xi_1,\xi_2,...,\xi_n\}\subset U(H)$ and $\lambda_1,\lambda_2,...,\lambda_n\in\mathbb{C}_{\qi}^+$ such that:
\begin{enumerate}
   \item $\mathcal{T}\xi_k=\lambda_k\xi_k$;
   \item  $\displaystyle \mathcal{T}u=\sum_{k=1}^n(u,\xi_k)\lambda_k\xi_k$ for any $u\in H$.
\end{enumerate}
\end{thm}

\begin{rem}
Theorem \ref{spectral-theorem-for-n-dimension} is a result for left eigenvalue of left operator in finite dimensional quaternionic Hilbert space. According to  Remark  \ref{leftORright}, this theorem is actually the same as Theorem 4.6 of  \cite{farenick2003spectral} by Farenick and Pidkowich.
%Indeed, $\lambda_k$ is not necessary to be restricted in $\mathbb{C}_{\qi}^+$, it can be any elements of $\theta(\lambda_k)$.
\end{rem}

\begin{thm}
 If $n$ is a positive integer and $H_0$ is an $n$-dimensional subspace of a quaternionic Hilbert space $H$, then
\begin{enumerate}
  \item  every one-to-one $\Q$-linear mapping of $\Q^n$ onto $H_0$ is a isomorphism, and
  \item  $H_0$ is closed, that is, $H_0$ is complete.
\end{enumerate}
\end{thm}

From Theorem \ref{spectal-theorem-positive-operator}, we see that every  eigenspace belonging to the non-zero eigenvalues of  a compact positive operator is finite-dimensional.  Noticing that every normal operator in finite-dimensional quaternionic Hilbert spaces can be diagonalized, we obtain the following spectral theorem for compact normal operator in quaternionic Hilbert spaces.

\begin{thm}\label{spectral-theorem-for-compact-normal}
Let $H$ be a (left) quaternionic Hilbert space and $\mathcal{T}\in \mathcal{B}(H)$ be a compact normal operator. Then there is a (possibly finite) sequence $\{\lambda_k\}$ of quaternions satisfying $|\lambda_1|\geq|\lambda_2|\geq\cdot\cdot\cdot$ and a sequence $\{\xi_k\} \in U(H)$
such that
\begin{enumerate}
  \item  $\{\lambda_k\}$ are eigenvalues of $\mathcal{T}$ and $|\lambda_k|\rightarrow 0$ if  $\{\lambda_k\}$ is infinite;
  %\item if $H_k=\left\{u\in H:Tu=pu,p\in \theta(\lambda_k) \right\}$, then  $\dim H_k<\infty$;
  %\item if $H_k=\left\{u\in H:Tu=pu,p\in \theta(\lambda_k) \right\}$,then $H_k\perp H_l$ for $k\neq l$;
  \item  $\mathcal{T}\xi_k=\lambda_k\xi_k$ and $\mathcal{K}\xi_k=|\lambda_k|^2\xi_k$ for every $k$, where $\mathcal{K}=\mathcal{T}\mathcal{T}^*$. For every $u\in H$, we have
  \begin{equation}\label{spectral expansion}
    \mathcal{T}u=\sum_{k=1}^\infty (u,\xi_k)\lambda_k\xi_k
  \end{equation}
and
  \begin{equation*}
\mathcal{K}u=\sum_{k=1}^\infty (u,\xi_k)|\lambda_k|^2\xi_k;
  \end{equation*}
  %\item if $P_k$ is the orthogonal projection from $H$ into $H_k$, then
  %\begin{equation*}
   % T=\sum_{k=1}^\infty \lambda_kP_k;
  %\end{equation*}
  \item  if $0$ is not a eigenvalue of $\mathcal{T}$, then $\{\xi_k\}$ is a orthonormal basis of $H$.
 % \item if $H_{\alpha_k}=\left\{u\in H:Tu=pu,p\in \theta(\alpha_k) \right\}$,then $\dim H_{\alpha_k}<\infty$ for $\alpha_k\neq0$.
\end{enumerate}
\end{thm}
\begin{proof}
Let $H_1,H_2,...$ be the eigenspaces belonging to the non-zero eigenvalues $\mu_1>\mu_2>...$ of the compact positive operator $\mathcal{K}=\mathcal{T}\mathcal{T}^*=\mathcal{T}^*\mathcal{T}$, and let $H_0=\mathrm{Ker}\mathcal{K}$. Since $\mathcal{K}\mathcal{T}=\mathcal{T}\mathcal{K}$, each $H_n(n\geq1)$ is invariant under $\mathcal{T}$, i.e. $\mathcal{T}(H_n)\subset H_n$. For $u\in H_0$ we have
\begin{equation*}
  \norm{\mathcal{T}u}^2=(\mathcal{T}u,\mathcal{T}u)=(\mathcal{T}^*\mathcal{T}u,u)=(\mathcal{K}x,x)=0
\end{equation*}
and so $\mathcal{T}|_{H_0}=0$. Furthermore, for $n\geq1$ the restriction of $\mathcal{T}$ to $H_n$ is normal; as $H_n$ is finite-dimensional, we know from Theorem \ref{spectral-theorem-for-n-dimension} that $H_n$ has an orthonormal basis consisting of eigenvectors of $\mathcal{T}$. Denote the dimension of $H_n(n\geq1)$ by $l_n$ and let $l_0=0,~s_n=\sum_{k=0}^nl_n$, then we can find $\lambda_{s_{n-1}+1},\lambda_{s_{n-1}+2},...,\lambda_{s_n}\in\Q$ and
$\{\xi_{s_{n-1}+1},\xi_{s_{n-1}+2},...,\xi_{s_n}\}\subset H_n$ such that
\begin{enumerate}
  \item  $\mathcal{T}\xi_k=\lambda_k\xi_k$ for every $s_{n-1}<k\leq s_n$;
  \item $\mathcal{K}\xi_k=|\lambda_k|^2\xi_k$ for every $s_{n-1}<k\leq s_n$;
  \item  $|\lambda_k|^2\in \theta(\mu_n)=\{\mu_n\}$, i.e. $|\lambda_k|^2=\mu_n$ for every $s_{n-1}<k\leq s_n$.
\end{enumerate}
 Taking the union of these bases, together with an orthonomal basis of $H_0$, we obtain an orthonormal basis consisting of eigenvectors of $\mathcal{T}$.
\end{proof}
\begin{rem}
Unlike the complex cases, $\lambda_k$ can be replaced any element in $\theta(\lambda_k)$ of Theorem \ref{spectral-theorem-for-compact-normal}.
\end{rem}

The authors  in \cite{Colombo2011Noncommutative} introduced the slice hyperholomorphic functional calculus. The  concept of  slice hyperholomorphic  opens a   new range of research for quaternion analysis. Moreover, the Fueter mapping theorem and  its inverse theorem reveal a deep relation among complex holomorphic, slice hyperholomorphic and monogenic functions. The notion of $S$-spectrum for   quaternionic operators, which arises naturally in the slice hyperholomorphic functional calculus, is a natural and elegant generalization of the classical notion of  spectrum and it provides a powerful tool for the study of spectral theorem for quaternionic operators. Based on the $S$-spectrum, the spectral theory of quaternionic operator has been developed  (see e.g. \cite{ghiloni2014spectral,Alpay2016unbounded,Alpay2016unitary,Colombo2016Schatten}).  Numerous celebrated spectral theorems in the classical case have been  carried out to the quaternionic case.

Based on $S$-spectrum, a spectral theorem for compact normal quaternionic operators was established in \cite{ghiloni2014spectral}. However, they did not sort the eigenvalues by norm in descending order. The spectral theorem for normal operators was generalized to  the unbounded operators by Alpay {\it et al.} \cite{Alpay2016unbounded}. Colombo {\it et al.} \cite{Colombo2016Schatten} presented  a singular value decomposition of compact quaternionic operator $T$ as:
 \begin{equation*}
    Tx=\sum_{n\in \mathbb{N}} \sigma_n\lambda_n \langle e_n,x \rangle,   ~~~~~\forall x\in H
 \end{equation*}
 where $\lambda_n >0$ are the singular values (not eigenvalues) of $T$, the vectors $\{e_n\}$ form an eigenbasis $|T|$, and $\sigma_n=We_n$ with $W$ unitary on $\ker W^{\bot}$ such that $T=W|T|$ (for detail, see Remark 3.4 of \cite{Colombo2016Schatten}). In the present paper, we obtain  a new and different spectral theorem, that is,
  \begin{equation*}
    \mathcal{T}u=\sum_{k=1}^\infty (u,\xi_k)\lambda_k\xi_k.
  \end{equation*}
Note that   $\{\lambda_k\}$ are eigenvalues of $\mathcal{T}$   and $\mathcal{T}\xi_k=\lambda_k\xi_k$.  Moreover, the eigenvalues are sorted by norm in descending order. This property will contribute to the settlement of   energy concentration problem in the following section.  Our method is   different from the $S$-spectrum based method and the newly obtained theorem is supplement  of the previous results.

%The proposed  method in present paper is useful to derive a spectral theorem for compact normal quaternionic operators even though   there is a question whether the method is  applicable to  other types of quaternionic operators.

It is   worth pointing out that  the authors in \cite{Colombo2016Schatten} also proposed the Schatten classes of quaternionic operators. Moreover, some characterizations of Scatten class operators were drawn. This topic is of great   importance in operator theory and it has huge potential in  physical and engineering applications.

We end up this part with an example to illustrate Theorem \ref{spectral-theorem-for-compact-normal}.

\begin{ex}
 Let $H=L^2([-\pi,\pi],\Q)$ be the quaternionic Hilbert space which consists of all square integrable $\Q$-valued functions defined on $[-\pi,\pi]$. Define the operator $T: F\mapsto f$ as follows:
 \begin{equation*}
   f(x)=T F(x)=\int_{-\pi}^{\pi}F(w)\overline{h(x,w)} dw=(F, h(x,\cdot))
 \end{equation*}
where $h(x,w)=-\qi \sin(wx)-\qj \cos(2wx)$. We conclude that
\begin{equation*}
     T^*F (x)  = \int_{-\pi}^{\pi}F(w) {h(x,w)} dw= -TF(x),
\end{equation*}
since
\begin{equation*}
  \begin{split}
    (F_1,\mathcal{T}^*F_2) =& \int_{-\pi}^{\pi} F_1(w) \overline{\left(\int_{-\pi}^{\pi} F_2(y)h(w,y)dy\right)}dw \\
           =& \int_{-\pi}^{\pi}\int_{-\pi}^{\pi} F_1(w)\overline{h(w,y)} ~\overline{F_2(y)}dw dy  \\
           =& \int_{-\pi}^{\pi} \mathcal{T}F_1(y)\overline{F_2(y)}dy= (\mathcal{T}F_1,F_2).
  \end{split}
\end{equation*}
That means $T^*=-T$. Therefore  $T^*T=T T^*$ which implies that $T$ is normal. Let $K=T^* T$, it follows that $K$ is   positive  and
\begin{equation*}
   KF(y)=\int_{-\pi}^{\pi}F(x)S(x,y)dx,
\end{equation*}
where $S(x,y)=\frac{\sin2\pi(x-y)-2\sin\pi(x-y)}{2(x-y)}+\frac{\sin2\pi(x+y)+2\sin\pi(x+y)}{2(x+y)}$.
As  functions of two variables, both $h$ and $S$ are continuous on $[-\pi,\pi]^2$. Therefore $T$ and $K$ are compact.
\end{ex}

%%%%%%%%%%%%%%%%%%%%%%%%%%%%%%%%%%%%%%%%%%%%%%%%%%%%%%%%%%%%%%%%
\section{Quaternionic Prolate Spheroidal  Wave Functions} \label{S4}

 In \cite{saitoh1983hilbert}, the author introduced the induced Hilbert spaces.  Motivated by this study,  we consider the integral transforms in quaternionic Hilbert spaces. The theory of   quaternionic reproducing kernel Hilbert spaces presented by  Alpay  {\it et al.}  \cite{Alpay2016SliceSchur} plays an important role in this part.

Let $\X$ be an arbitrary set and $\mathfrak{F}(\X)$ a $\Q$-linear space composed of all $\Q$-valued signals on $\X$. Let $H$ be a (possibly finite-dimensional) quaternionic Hilbert space with inner product $(\cdot,\cdot)_H$, and $E:\X\rightarrow H$ be a vector-valued function from $\X$ into $H$. Consider the left $\Q$-linear mapping $\mathcal{T}$ from $H$ into $\mathfrak{F}(\X)$ defined by
\begin{equation*}
  f(\bx)=(\mathcal{T}F)(\bx)=(F,E(\bx))_H,
\end{equation*}
where $\mathcal{T}F=f$, $F\in H$, $f\in \mathfrak{F}(\X)$.

Let $\mathcal{H}$ and $N$ denote the range and null space of $\mathcal{T}$. Then $N$ is a closed subspace of $H$. Denote $M=N^\perp$, then $H=M\oplus N$.

\begin{thm}\label{reproducing-kernel-Hilbert-spaces}
$(\mathcal{H},(\cdot,\cdot)_{\mathcal{H}})$ is a quaternionic Hilbert space that is isometric to $(M,(\cdot,\cdot)_H)$, where
\begin{equation}\label{inner-product-of-range}
  (f,g)_{\mathcal{H}}=(\mathcal{T}F,\mathcal{T}G)_{\mathcal{H}}=(P_MF,P_MG)_H.
\end{equation}
Moreover, $(\mathcal{H},(\cdot,\cdot)_{\mathcal{H}})$ is a quaternion reproducing kernel Hilbert space with the kernel $S(\bx,\y)$ defined as
\begin{equation*}
  S(\bx,\y)=(E(\y),E(\bx))_H.
\end{equation*}
\end{thm}

\begin{proof}
It is easy to see that $\mathcal{T}$ is a bijection from $M$ into $H$. It implies that for every $f\in \mathcal{H}$, there exists a unique element $F^M_f\in M$ such that $\mathcal{T}F^M_f=f$. We can show that $\|F^M_f\|^2_{H}=\inf\{\|F\|^2_H:f=\mathcal{T}F\}$. It is obviously that $\|F^M_f\|^2_{H}\geq \inf\{\|F\|^2_H:f=\mathcal{T}F\}$. On the other hand, for every $F\in\{F\in H:\mathcal{T}F=f\}$, it can be decomposed as $F=F^M_f+F'$,where $F'\in N$. Then
\begin{equation*}
  \begin{split}
   \|F\|^2_H  =  (F^M_f+F',F^M_f+F')_H
   =& \|F^M_f\|^2_{H}+(F',F^M_f)_H+(F^M_f,F')_H+\|F'\|^2_H\\
   =& \|F^M_f\|^2_{H}+\|F'\|^2_H
   \geq  \|F^M_f\|^2_{H}.
  \end{split}
\end{equation*}
Thus $\inf\{\|F\|^2_H:f=\mathcal{T}F\} \geq \|F^M_f\|^2_{H}$. Hence $\|f\|^2_{\mathcal{H}}=\|F^M_f\|^2_{H}=\inf\{\|F\|^2_H:f=\mathcal{T}F\}$ and $\mathcal{T}\mid_M$ is an isometry between $(M,(\cdot,\cdot)_H)$ and $(\mathcal{H},(\cdot,\cdot)_{\mathcal{H}})$.

Next, we prove that for every $\bx\in \X, E(\bx)\in M$. Noticing that when $F\in N$, then for every $\bx\in \X$ we have
\begin{equation*}
f(\bx)=(F,E(\bx))_H=0.
\end{equation*}
That means for every $\bx\in \X$, $E(\bx)\in N^\perp=M$ and $P_ME(\bx)=E(\bx)$.
Since $S(\bx,\y)=(E(\y),E(\bx))_H=\mathcal{T}(E(\y))(\bx)$, thus
\begin{equation*}
  \begin{split}
  (f,S(\cdot,\y))_{\mathcal{H}}=&(\mathcal{T}F,\mathcal{T}(E(\y)))_{\mathcal{H}}=(P_MF,P_ME(\y))_H  \\
  =&(F^M_f,E(\y))_H=\mathcal{T}F^M_f(\y)=f(\y)
  \end{split}
\end{equation*}
which completes the proof.
\end{proof}

For $\Q$-valued signals $f,g:\D\rightarrow \Q$ where $\D$ is a compact connected subset of $\mathbb{R}^d$, we can define the $\Q$-valued inner product
$ (f,g)=\int_\D f(\w)\overline{g(\w)}d\w.$
The left  $\Q$-linear quaternionic Hilbert space $L^2(\D,\Q)$ consists of all $\Q$-valued signals which are square-integrable on $\D$:
$$L^2(\D,\Q)=\left\{f|f:\D\rightarrow \Q,\|f\|:=\int_\D |f(\w)|^2d\w<\infty\right\}.$$
Now we apply Theorem \ref{reproducing-kernel-Hilbert-spaces} to a specific case. Let $H=L^2(\D,\Q)$ and $\X$ be an open, connected subset of $\mathbb{R}^d$ containing $\D$. The  $\Q$-valued function $E(\w,\bx)$ on $\D\times \X$ satisfying  $E(\w,\bx)\in L^2(\D,\Q)$ for any $\bx\in \X$. In the next, we consider the integral transform of $F\in L^2(\D,\Q)$,
\begin{equation}\label{integral-transform}
  f(\bx)=(\mathcal{T}F)(\bx)=\int_\D F(\w)\overline{E(\w,\bx)}d\w.
\end{equation}
\begin{thm}\label{compact}
  Suppose that $E(\w,\bx)$ is square-integrable on $\D\times \D$. Then $\mathcal{T}\in \mathcal{B}_0(H)$  where $H=L^2(\D,\Q)$ and $\mathcal{T}$ is given by (\ref{integral-transform}).
\end{thm}
\begin{proof}
Clearly, $\mathcal{T}$ is linear. Since
  \begin{equation*}
    \begin{split}
   \int_\D|f(\bx)|^2 d\bx =& \int_\D \left|\int_\D F(\w)\overline{E(\w,\bx)}d\w\right|^2d\bx \\
     \leq& \int_\D\left( \int_\D\left|F(\w)\right|^2d\w\int_\D\left|E(\w,\bx)\right|^2d\w\right)d\bx \\
     \leq&  \int_\D\left|F(\w)\right|^2d\w\int_{\D^2}\left|E(\w,\bx)\right|^2d\w d\bx<\infty.
    \end{split}
  \end{equation*}
  Therefore $f(\bx)\in H$.  Since $\norm{\mathcal{T}F}^2\leq \alpha^2\norm{F}^2$, where
           $$\alpha=\left(\int_{\D^2}\left|E(\w,\bx)\right|^2d\w d\bx\right)^{\frac{1}{2}},$$
      then     $\mathcal{T}$ is bounded.

 Suppose $\{e_n(\w)\}$ is an orthonormal basis of $H$. % $H$ is separable
  It is easily checked that ${\overline{e_n(\bx)}e_m(\w)}$ is an orthonormal basis of $L^2(\D^2,\Q)$.
  Thus
   $ E(\w,\bx)=\sum_{m,n} c_{mn} {\overline{e_n(\bx)}e_m(\w)}$
  where
  \begin{equation*}
 c_{mn}=\int_{\D^2} E(\w,\bx)\overline{\overline{e_n(\bx)}e_m(\w)}d\w d\bx \\
    =\int_{\D^2} E(\w,\bx)\overline{e_m(\w)}e_n(\bx)d\w d\bx.
 \end{equation*}
  From Parseval's identity, we have $\sum\limits_{m,n}|c_{mn}|^2=\alpha^2<\infty$.
  For each $e_m(\w)\in H$, $g_m(\bx)=\overline{(\mathcal{T}e_m)(\bx)}\in H$, then
  \begin{equation*}
    g_m(\bx)=\sum_n d_{mn}e_n(\bx),
  \end{equation*}
  where
  \begin{equation*}
    \begin{split}
       d_{mn} =& \int_\D g_m(\bx)\overline{e_n(\bx)}d\bx
     = \int_\D \left(\overline{\int_\D e_m(\w)\overline{E(\w,\bx)}d\w}\right)\overline{e_n(\bx)}d\bx  \\
     =& \int_{\D^2} E(\w,\bx)\overline{e_m(\w)}e_n(\bx)d\w d\bx = c_{mn}.
    \end{split}
  \end{equation*}
  Since $\sum\limits_{m,n}|c_{mn}|^2$ converges, then
  \begin{equation*}
    \lim_{m\rightarrow\infty}\norm{\mathcal{T}e_m}^2=\lim_{m\rightarrow\infty}\norm{g_m}^2=\lim_{m\rightarrow\infty}\sum_{n=1}^{\infty}|c_{mn}|^2=0.
  \end{equation*}
   In fact, for the same reason, we have
  \begin{equation*}
    \lim_{k\rightarrow\infty}\sum_{m>k}\norm{\mathcal{T}e_m}^2=\lim_{k\rightarrow\infty}\sum_{m>k}\sum_{n=1}^{\infty}|c_{mn}|^2=0.
  \end{equation*}
  For every $\bx\in \D$, $E(\w,\bx) \in H$, thus
  \begin{equation*}
    E(\w,\bx)=\sum_{n=1}^\infty(E(\cdot,\bx),e_n(\cdot))e_n(\w)=\sum_{n=1}^\infty\overline{\mathcal{T}e_n(\bx)}e_n(\w).
  \end{equation*}

  Truncate the kernel $E(\w,\bx)$  by
  \begin{equation*}
     E_k(\w,\bx)=\sum_{n=1}^k\overline{\mathcal{T}e_n(\bx)}e_n(\w)=\sum_{n=1}^kg_n(\bx)e_n(\w),
  \end{equation*}
These give the finite-rank operators
%A continuous linear operator is of finite rank if its image is finite-dimensional. Note that a finite-rank
%operator is compact, since all balls are pre-compact in a finite-dimensional Banach space.
\begin{equation*}
  (\mathcal{T}_kF)(\bx)=\int_\D F(\w) \overline{E_k(\w,\bx)}d\w,~~~\mathrm{for ~every} ~~~F(\w)\in H.
\end{equation*}
We claim that $\{\mathcal{T}_k\}$ tend to $\mathcal{T}$ in operator norm. Indeed, for every $F(\w)=\sum_n b_n e_n(\w)$ in $H$, we have
\begin{equation*}
  \begin{split}
 ((\mathcal{T}-\mathcal{T}_k)F)(\bx)=&(F,E(\cdot,\bx))-(F,E_k(\cdot,\bx)) \\
   =& \sum_{n=1}^{\infty} (F,e_n)(e_n,E(\cdot,\bx))-(F,\sum_{n=1}^{k}\overline{\mathcal{T}e_n(\bx)}e_n) \\
   =& \sum_{n=1}^{\infty} (F,e_n)\mathcal{T}e_n(\bx)-\sum_{n=1}^{k}(F,e_n)\mathcal{T}e_n(\bx) \\
   =& \sum_{n>k} b_n \mathcal{T}e_n(\bx).
  \end{split}
\end{equation*}
Therefore
\begin{equation*}
  \begin{split}
  \norm{(\mathcal{T}-\mathcal{T}_k)F} =& \norm{\sum_{n>k} b_n \mathcal{T}e_n}
   \leq  \sum_{n>k}|b_n| \norm{\mathcal{T}e_n} \\
   \leq& \left(\sum_{n>k}|b_n|^2\right)^{\frac{1}{2}}\left(\sum_{n>k}\norm{\mathcal{T}e_n}^2\right)^{\frac{1}{2}}
   \leq  \norm{F} \left(\sum_{n>k}\norm{\mathcal{T}e_n}^2\right)^{\frac{1}{2}}.
  \end{split}
\end{equation*}
As $\lim_{k\rightarrow\infty}\sum_{n>k}\norm{\mathcal{T}e_n}^2=0$, we have
$\lim_{k\rightarrow\infty}\mathcal{T}_k=\mathcal{T}$.

Since every finite-rank operator is compact, and $\mathcal{B}_0(H)$ is a closed subset of $ \mathcal{B}(H)$, thus $\mathcal{T}$ is compact.
\end{proof}

\begin{rem}
By Theorem \ref{reproducing-kernel-Hilbert-spaces} and Theorem \ref{compact}, we see that $\mathcal{T}$ can be regarded an operator  of both  ${\mathcal{B}}(H,\mathcal{H})$ and $\mathcal{B}(H)$.
\end{rem}

%%%%%%%%%%%%%%%%%%%%%%%%%%%%%%%%%%%%%%%%%%%%%%%%%%%%%%%%%%%%%%%%
\begin{thm}\label{main-results-GQPSWFs}
%Let $\D$ be a compact connected subset of $\mathbb{R}^N$ with $m(\D)>0$, and $\X$ be an open, connected subset of $\mathbb{R}^N$ containing $\D$, i.e. $\D\subset \X\subset \mathbb{R}^N$.
  Let $E(\w,\bx)$ be a $\Q$-valued continuous function defined on $\X^2\subset\mathbb{R}^{2d}$ satisfying the following conditions:
\begin{enumerate}
  \item $E(\w,\bx)=E(\bx,\w)$ for every $(\w,\bx)\in \X^2$.
  \item  $A=\{E(\w,\bx)\}_{\bx\in \X}\subset H=L^2(\D,\Q)$ and $A^{\perp}=\{0\}$.
  \item  $S(\bx,\y)= \int_\D E(\w,\y) \overline{E(\w,\bx)} d\w\in \mathbb{R} $ for every $(\bx,\y)\in \X^2$;
  \item  ${\int_\D E(\w,\y) \overline{E(\w,\bx)} d\w= \int_\D \overline{E(\w,\bx)} E(\w,\y)d\w}$  for every $(\bx,\y)\in \X^2$.
\end{enumerate}
Then we can find a countably infinite set of $\Q$-valued signals ${\{\phi_n(\bx)\}}_{n=1}^{\infty}$ called quaternionic prolate spheroidal  wave functions(QPSWFs) and a set of $\Q$-valued numbers
$|\lambda_1|\geq|\lambda_2|\geq\cdot\cdot\cdot$ with the following properties:

\begin{enumerate}

  \item  The ${ \phi_n(\bx) }~(n=1,2,\cdot\cdot\cdot)$  are $\D$-bandlimited in the sense of transformation $\mathcal{T}$,  orthonormal  and complete in $\mathcal{H}$, where $\mathcal{H}$ is a quaternionic Hilbert space  with the inner product defined by (\ref{inner-product-of-range}):
\begin{equation*}
(\phi_m,\phi_n)_{\mathcal{H}}=\delta_{mn}.
\end{equation*}
\item
The ${ \phi_n(\bx) }~(n=1,2,\cdot\cdot\cdot)$ are orthogonal and complete in $L^2(\D,\Q)$:
 \begin{equation*}
  \int_\D \phi_m (\w) \overline{\phi_n(\w)}  d\w=\mu_n\delta_{mn},~~\mathrm{with}~ ~\mu_n=|\lambda_n|^2.
\end{equation*}
%where $\{\lambda_n\}$ are the eigenvalues of a compact operator  $L\in \mathcal{B}_0(H)$.\\
%   ${\{\varphi_n\}}_{n=1}^{\infty}$   are solutions of the following two integral equations,

\item  For every $\bx\in \X$, we have
\begin{align}
    &\int_\D \phi_n(\w)\overline{E(\w,\bx)}d\w=\lambda_n\phi_n(\bx),  \label{aa} \\
   &  \int_\D\phi_n (\y) S(\bx,\y) \,d\y=\mu_n\phi_n (\bx). \label{bbb}
\end{align}
\item
 The ${ \phi_n(\bx) }~(n=1,2,\cdot\cdot\cdot)$ are     uniformly continuous on $\X_1$, where $\X_1$ being any compact subset of $\X$; and, hence $\phi_n(\bx)$ is continuous on $\X$.
  \item  $E(\w,\bx)$ and $S(\bx,\y)$ can be expanded by $\phi_n$:
  \begin{align}
    &  E(\w,\bx)=\sum_{n=1}^{\infty}\overline{\phi_n(\bx)}\lambda_n^{-1}\phi_n(\w)  ~~~(\w,\bx)\in \D\times \X,  \label{expansion-of-h}\\
   &  S(\bx,\y)=\sum_{n=1}^{\infty}\overline{\phi_n(\y)} \phi_n(\bx)~~~(\bx,\y)\in \X^2. \label{expansion-of-K}
\end{align}
For any fixed $\bx \in \X$, the series (\ref{expansion-of-h}) converges in the norm.
The series (\ref{expansion-of-K}) converges   absolutely on $\X\times \X$. Furthermore, if $E(\w,\bx)$ is bounded on $\D\times \X$, then the series (\ref{expansion-of-K}) converges uniformly on $\X_1\times \X_2$, where $\X_1,\X_2$ being any compact subsets of $\X$.

\item  If there exists a sequence of points
   ${\{\bx_{n}\}}{\subset}\X $ such that $\{E(\w,\bx_{n})\}$ is
 an orthonormal basis of $L^2 (\D,\Q)$, then for any $f\in  \mathcal{H}$,
\begin{align}
    &  f(\bx)=\sum_{n=1}^\infty f(\bx_{n})S(\bx,\bx_{n}),  \label{shannon sampling} \\
   &  f(\bx)=\sum_{n=1}^\infty \left(\sum_{m=1}^\infty f(\bx_m)\overline{\phi_n(\bx_m)}\right) \phi_n(\bx). \label{sampling with QPSWFs}
\end{align}
\item   The signals $\{\phi_n\}$  satisfy the discrete orthogonality relation:
\begin{align*}
    &    \sum_{n=1}^{\infty}\overline{\phi_n(\bx_l)} \phi_n(\bx_m)=\delta_{ml}; \\
   &    \int_\D S(\bx,\bx)d\bx=\sum_{n=1}^{\infty}\mu_n.
\end{align*}
  \item  For every $f(\bx)\in \mathcal{H}$ with $\norm{f}_\mathcal{H}\neq 0$, we form the ratio
\begin{equation}
  \displaystyle \beta_f= \frac{\norm{f(\bx)}_H^2}{\norm{f(\bx)}_\mathcal{H}^2}.
\end{equation}
If we put $\widetilde{\beta}=\sup\{\beta_f:f\in\mathcal{H},\norm{f}_\mathcal{H}\neq0\}$, then $\widetilde{\beta}=\mu_1$.
\end{enumerate}
\end{thm}
\begin{proof}
%{\em\romannumeral1.}
Since ${\overline{S(\bx,\y)}=\int_\D E(\w,\bx) \overline{E(\w,\y)} d\w=S(\y,\bx)}$, and $S$ is   real, then  $S(\bx,\y)=S(\y,\bx)$.

 We consider the integral transform $\mathcal{T}$ defined by (\ref{integral-transform}). Since $E(\w,\bx)$ is continuous on compact set $\D\times \D$, by Theorem \ref{compact}, we have $\mathcal{T}\in \mathcal{B}_0(H)$ when $\bx$ is restricted to $\D$. We now prove that $\mathcal{T}$ is normal.
Consider its adjoint operator $\mathcal{T}^*$. Indeed,
\begin{equation*}
  \mathcal{T}^*F(\bx)=\int_\D F(\w)E(\bx,\w)d\w,
\end{equation*}
since
\begin{equation*}
  \begin{split}
    (F_1,\mathcal{T}^*F_2)_H =& \int_\D F_1(\w) \overline{\left(\int_\D F_2(\y)E(\w,\y)d\y\right)}d\w \\
           =& \int_{\D^2} F_1(\w)\overline{E(\w,\y)} ~\overline{F_2(\y)}d\w d\y  \\
           =&  \int_\D \mathcal{T}F_1(\y)\overline{F_2(\y)}d\y= (\mathcal{T}F_1,F_2)_H.
  \end{split}
\end{equation*}
Therefore
\begin{equation*}
  \begin{split}
     \mathcal{T}(\mathcal{T}^*F)(\bx) =& \mathcal{T}\left(\int_\D F(\w)E(\y,\w)d\w\right)(\bx) \\
     =& \int_{\D^2}F(\w)E(\y,\w)\overline{E(\y,\bx)}d\y d\w\\
     =&  \int_\D F(\w)S(\bx,\w)d\w,
  \end{split}
\end{equation*}
and
\begin{equation*}
  \begin{split}
  \mathcal{T}^*(\mathcal{T}F)(\bx) =&\mathcal{ T}^*\left(\int_\D F(\w)\overline{E(\w,\y)}d\w\right)(\bx) \\
     =& \int_{\D^2}F(\w)\overline{E(\w,\y)}E(\bx,\w)d\w d\y\\
     =&  \int_\D F(\y)S(\bx,\y)d\y.
  \end{split}
\end{equation*}
Hence, $\mathcal{T}^*\mathcal{T}=\mathcal{T}\mathcal{T}^*$. Let $\mathcal{T}^*\mathcal{T}=\mathcal{T}\mathcal{T}^*=\mathcal{K} $, then
\begin{equation}\label{positiveK}
g(\bx)= \mathcal{K}F(\bx)=\int_D F(\y)S(\bx,\y)d\y,~~~~\bx\in \D.
\end{equation}
 Thus, by Theorem {\ref{spectral-theorem-for-compact-normal}}, there is a sequence $\{\lambda_n\}$ of quaternions satisfying $|\lambda_1|\geq|\lambda_2|\geq\cdot\cdot\cdot$ and a sequence $\{\Phi_n\}_{n=1}^\infty \in U(H)$,
such that
\begin{equation}\label{eigenfunction1}
  \int_\D \Phi_n(\w)\overline{E(\w,\bx)}d\w=\lambda_n\Phi_n(\bx),~\bx\in \D
\end{equation}
and
\begin{eqnarray} \label{eigenfunction2}
 { \int_\D \Phi_n(\y)S(\bx,\y)d\y=|\lambda_n|^2\Phi_n(\bx)=\mu_n\Phi_n(\bx),~\bx\in \D.}
\end{eqnarray}

Since $A^{\perp}=0$, which implies the  null space of of $\mathcal{T}$ is $\{0\}$ and $\lambda_n \neq0$ for every $n$. Hence, $\{\Phi_n\}_{n=1}^\infty$ is an orthonormal basis of $H$.

From Theorem \ref{reproducing-kernel-Hilbert-spaces}, we see that $\mathcal{T}$ is an isometry between $(H,(\cdot,\cdot)_H)$  and $(\mathcal{H},(\cdot,\cdot)_{\mathcal{H}})$. Set
%{\lambda_n^{-1}}
\begin{equation*}
  \phi_n(\bx)=\left(\mathcal{T}\Phi_n\right)(\bx),~~~\bx\in \X.
\end{equation*}

Then $\{\phi_n\}_{n=1}^{\infty}$ is an orthonormal basis of $\mathcal{H}$ as $\{\Phi_n\}_{n=1}^\infty$ is an orthonormal basis of $H$.
Moreover,   $\phi_n(\bx)=\lambda_n\Phi_n(\bx)$ for every $\bx\in \D$ from (\ref{eigenfunction1}). Thus  $\phi_n(\bx)$ is an extension of $\lambda_n\Phi_n$. Furthermore,
\begin{equation*}
  (\phi_m,\phi_n)_{\mathcal{H}}=(\Phi_m,\Phi_n)_H=\delta_{mn},
\end{equation*}

and
\begin{equation*}
      \int_\D \phi_m (\w) \overline{\phi_n(\w)}  d\w= \int_\D\lambda_m \Phi_m (\w) \overline{\lambda_n\Phi_n(\w)} d\w
  =\lambda_m(\Phi_m,\Phi_n)_H\overline{\lambda_n}=\mu_n\delta_{mn}.
\end{equation*}

 If we regard $\mathcal{T}$ as an operator of $B(H,\mathcal{H})$, then $\mathcal{K}=\mathcal{T}\mathcal{T}^*$ is also an operator of $B(H,\mathcal{H})$. Thus $\bx\in \D$ of (\ref{positiveK}) can be replaced by $\bx\in \X$.  Since $\mathcal{T}^*\Phi_n=\overline{\lambda}_n\Phi_n$, then $$\overline{\lambda}_n\mathcal{T}\Phi_n=\mathcal{T}(\overline{\lambda}_n\Phi_n)=\mathcal{T}\mathcal{T}^*\Phi_n=\mathcal{K}\Phi_n.$$ Hence, $\phi_n(\bx)=\overline{\lambda}_n^{-1}\mathcal{K}\Phi_n(\bx)$.
From the definition of $\mathcal{T}$ and $\mathcal{K}$, we easily find
\begin{equation*}
  \begin{split}
   \lambda_n\phi_n(\bx)=\lambda_n\mathcal{T}\Phi_n(\bx)= \lambda_n\int_\D \Phi_n(\w)\overline{E(\w,\bx)}d\w
  =&\int_\D \lambda_n \Phi_n(\w)\overline{E(\w,\bx)}d\w\\
  =&\int_\D \phi_n(\w)\overline{E(\w,\bx)}d\w
  \end{split}
\end{equation*}
and
\begin{equation*}
  \mu_n\phi_n(\bx) = \lambda_n \overline{\lambda}_n\overline{\lambda}_n^{-1}K\Phi_n(\bx)
   =   \int_\D \lambda_n \Phi_n(\y)S(\bx,\y)d\y
   =\int_\D  \phi_n(\y)S(\bx,\y)d\y.
\end{equation*}
Then we obtain (\ref{aa}) and (\ref{bbb}).

For given $F\in H$, $f(\bx)=\mathcal{T}F(\bx)$, we have
\begin{equation*}
  \begin{split}
  |f(\bx)-f(\y)| =& \left| \int_\D F(\w)\left[ \overline{E(\w,\bx)-E(\w,\y)}\right]d\w\right| \\
   \leq&  \norm{F}_H \left[\int_\D |E(\w,\bx)-E(\w,\y)|^2d\w\right]^{\frac{1}{2}}
  \end{split}
\end{equation*}
Since $E(\w,\bx)$ is continuous on compact set $\D\times \X_1$, then $E(\w,\bx)$ is uniformly continuous on $\D\times \X_1$.
Thus  $\forall \varepsilon>0$, $\exists\delta>0$, when $\bx,\y\in X_1$ and $d(\bx,\y)<\delta$,
\begin{equation*}
 |E(\w,\bx)-E(\w,\y)|<\frac{\varepsilon}{\norm{F}_H\cdot (\varrho(\D))^{\frac{1}{2}}}
\end{equation*}
holds for every $\w\in \D$, where $\varrho(\D)$ is the Lebesgue measure of $\D$.
Therefore
\begin{equation*}
|f(\bx)-f(\y)|<\varepsilon.
 \end{equation*}
Hence, $f(\bx)$ is uniformly continuous on $\X_1$.

%{\em\romannumeral2.}
Since $\{\Phi_n\}_{n=1}^\infty$ is an orthonormal basis of $H$, thus
\begin{equation*}
  \begin{split}
  E(\w,\bx) {=} \sum_{n=1}^{\infty}(E(\cdot,\bx),\Phi_n(\cdot))_H\Phi_n(\w)
  =& \sum_{n=1}^{\infty}\overline{(\Phi_n(\cdot),E(\cdot,\bx))_H}\Phi_n(\w)   \\
     =& \sum_{n=1}^{\infty}\overline{\phi_n(\bx)}\lambda_n^{-1}\phi_n(\w),
  \end{split}
\end{equation*}
and
\begin{equation*}
  \begin{split}
    S(\bx,\y) =  (E(\cdot,\y),E(\cdot,\bx))_H
   =& \sum_{n=1}^{\infty}(E(\cdot,\y),\Phi_n(\cdot))_H(\Phi_n(\cdot),E(\cdot,\bx))_H \\
   =&  \sum_{n=1}^{\infty}\overline{\phi_n(\y)} \phi_n(\bx).
  \end{split}
\end{equation*}
If $E(\w,\bx)$ is continuous and bounded on $\D\times \X$, suppose that $|E(\w,\bx)|<M_1$ for every $(\w,\bx)\in \D\times \X$. Then
\begin{eqnarray*}
  0<S(\bx,\bx) %&=& \sum_{n=1}^{\infty}\overline{\Phi_n(p)} \Phi_n(p) \\
         &=&  \sum_{n=1}^{\infty}|\phi_n(\bx)|^2
         = (E(\cdot,\bx),E(\cdot,\bx))_H \\
         &=& \int_\D|E(\w,\bx)|^2 d\w<M_1^2\cdot \varrho(\D)<\infty.
\end{eqnarray*}
For any fixed $\y\in \X_1$ and given positive $\varepsilon$, we have
\begin{equation*}
\begin{split}
  \left(\sum_{n=k}^{m}\left|\overline{\phi_n(\y)} \phi_n(\bx)\right| \right)^2 \leq& \sum_{n=k}^{m}\left| \phi_n(\y) \right|^2\sum_{n=k}^{m}\left| \phi_n(\bx) \right|^2
   \leq  \sum_{n=k}^{m}\left| \phi_n(\y) \right|^2 S(\bx,\bx)  \\
   \leq& \varepsilon M_1^2\cdot \varrho(\D)~~\mathrm{for}~ k,m\geq N(\y,\varepsilon).
\end{split}
\end{equation*}
As $E(\w,\bx)$ is continuous, it is  easily seen that $S(\bx,\y)$ is continuous. Therefore, by Dini's theorem, the convergence of  $\sum_{n=1}^{\infty}|\phi_n(\y)|^2$ to  $S(\y,\y)$   on $\X_2$ must be uniform. The dependence upon $\y$ of $N(\y,\varepsilon)$ is thus actually extrinsic, from which fact the desired uniform (and absolute) convergence of the series to $S(\bx,\y)$ on $\X_1\times \X_2$ immediately follows.

%{\em\romannumeral3.}
Since $\left\{E_n(\w)=E(\w,\bx_n)\right\}$ is an orthonomal basis of $H$, we have
$F(\w)=\sum_{n} (F,E_n)_H E_n(\w)$.
Set $F_m(\w)=\sum_{n\leq m} (F,E_n)_H E_n(\w)$
then $F_m$ converges to $F$ in the norm as $m \rightarrow \infty$. Note that $f(\bx_n)=(F,E_n)_H$ and $S(\bx,\bx_n)=(E_n,E(\cdot,\bx)_H$, therefore
\begin{equation*}
\mathcal{T}F_m(\bx)=\sum_{n\leq m}(F,E_n)_H (E_n,E(\cdot,\bx)_H
 =\sum_{n\leq m} f(\bx_n)S(\bx,\bx_n).
\end{equation*}
We conclude that $\mathcal{T}F_m(\bx)$ converges uniformly to $\mathcal{T}F(\bx)$ on $\X$, that is, (\ref{shannon sampling}) converges uniformly on $\X$, since
\begin{equation*}
  \begin{split}
  \left|\mathcal{T}F(\bx)-\mathcal{T}F_m(\bx)\right|^2 =& \left|\int_\D [F(\w)-F_m(\w)]\overline{E(\w,\bx)}d\w\right|^2 \\
    \leq&  \int_\D \left|F(\w)-F_m(\w)\right|^2d\w \cdot \int_\D  \left|E(\w,\bx)\right|^2d\w \\
    \leq& \norm{F-F_m}_H^2 \cdot M_1^2\cdot \varrho(\D).
  \end{split}
\end{equation*}

Since $\left\{E_n(\w)\right\}$ is an orthonomal basis of $H$, then
\begin{equation*}
  \begin{split}
   f(\bx) =&    (F,E_\bx)_H
       =   \sum_{n=1}^\infty (F,\Phi_n)_H(\Phi_n,E_\bx)_H\\
       =&   \sum_{n=1}^\infty (F,\Phi_n)_H \phi_n(\bx)
       =    \sum_{n=1}^\infty \left(\sum_{m=1}^\infty(F,E_m)_H(E_m,\Phi_n)_H\right) \phi_n(\bx) \\
       =&   \sum_{n=1}^\infty \left(\sum_{m=1}^\infty f(\bx_m)\overline{\phi_n(\bx_m)}\right)\phi_n(\bx).
  \end{split}
\end{equation*}
%{\em\romannumeral4.}
On the one hand,
$ S(\bx_m,\bx_l) = \sum_{n=1}^{\infty}\overline{\phi_n(\bx_l)} \phi_n(\bx_m)$;
on the other hand,
$S(\bx_m,\bx_l) = (E_l,E_m)_H=\delta_{ml}$.
Therefore,
 \begin{equation*}
   \sum_{n=1}^{\infty}\overline{\phi_n(\bx_l)} \phi_n(\bx_m)=\delta_{ml} .
 \end{equation*}

 Since
$S(\bx,\bx)=\sum_{n=1}^{\infty}|\phi_n(\bx)|^2$
converges uniformly on $\D$, thus
\begin{equation*}
    \int_\D S(\bx,\bx)d\bx=\int_\D \sum_{n=1}^{\infty}|\phi_\n(\bx)|^2 d\bx
    =\sum_{n=1}^{\infty}\int_\D |\phi_n(\bx)|^2 d\bx=\sum_{n=1}^{\infty}\mu_n.
\end{equation*}

%{\em\romannumeral5.}
 Suppose that $f=\mathcal{T}F$. From (\ref{inner-product-of-range}), we have
 \begin{equation*}
   \norm{f}_\mathcal{H}^2=\norm{F}_H^2=\sum_{n=1}^\infty \left|a_n\right|^2,
 %  \sum_{n=1}^\infty \left|(F,\varphi_n)\right|^2=
 \end{equation*}
 where $a_n=(F,\Phi_n)_H$.
As $f(\bx)\in L^2 (\D,\Q)$ when $\bx$ is restricted in $\D$, it can be expanded into a series
$f(\bx)=\sum_{n=1}^{\infty}b_n\Phi_n(\bx)$
 where
\begin{equation*}
  \begin{split}
   b_n =& \int_\D f(\bx)\overline{\Phi_n(\bx)}d\bx
    = \int_\D \left(\int_\D F(\w) \overline{E(\w,\bx)}d\w \right) \overline{\Phi_n(\bx)}d\bx\\
    =&  \int_\D F(\w) \left( \overline{\int_\D \Phi_n(\bx)E(\w,\bx)d\bx }\right) d\w
    =   \int_\D F(\w) \left( \overline{T^* \Phi_n(\w)}\right) d\w \\
    %&=&  \int_T F(t) \left( \overline{\overline{\lambda_n} \varphi_n(t)}\right) dt \\
    =&  \left( \int_\D F(\w) \overline{ \Phi_n(\w)} d\w\right)\lambda_n =a_n\lambda_n.
  \end{split}
\end{equation*}
 Thus
 \begin{equation*}
   \norm{f}_H^2=\sum_{n=1}^\infty \left|a_n\right|^2\left|\lambda_n\right|^2=\sum_{n=1}^\infty \left|a_n\right|^2\mu_n.
 \end{equation*}
Since  $\{\mu_n\} $ is monotonically decreasing and positive, Therefore
\begin{equation*}
   \beta_f =\frac{ \sum_{n=1}^\infty \left|a_n\right|^2\mu_n}{\sum_{n=1}^\infty \left|a_n\right|^2}
   %= \frac{ \sum_{n=n_1}^\infty \left|a_n\right|^2\mu_n}{\sum_{n=n_1}^\infty \left|a_n\right|^2}
    %&\leq&\frac{ \sum_{n=n_1}^\infty \left|a_n\right|^2\mu_{n_1}}{\sum_{n=n_1}^\infty \left|a_n\right|^2}  \\
 \leq   \mu_1 .
\end{equation*}
 Taking $f(\bx)=\phi_1(\bx)$, we have $\beta_{\phi_1}= \mu_1$. Thus $\widetilde{\beta}= \mu_1$.
\end{proof}

\begin{ex}
Let $\tau, \sigma>0$, $\D=[-\tau,\tau]^2$, $\X=\mathbb{R}^2$, $E(\w,\bx)=\frac{1}{2\tau}\e^{-\qi \sigma x_1\omega_1/ \tau}\e^{-\qj \sigma x_2\omega_2/ \tau}$.
 Then $$S(\bx,\y)=\frac{\sin \sigma (x_1-y_1)\sin \sigma (x_2-y_2)}{ \sigma(x_1-y_1) \sigma(x_2-y_2)}.$$ We consider the following finite modified inverse quaternion Fourier transform
 \begin{equation*}
 f(x_1,x_2)= \mathcal{T}F(\bx)
 =\frac{1}{2\tau}\int_{\D}F(\omega_1,\omega_2)\e^{\qj \sigma x_2\omega_2/ \tau}\e^{\qi \sigma x_1\omega_1/ \tau}d\omega_1d\omega_2
 \end{equation*}
where $F(\omega_1,\omega_2)\in L^2(\D,\Q)$. It is easy to see that $f(x_1,x_2)$ is $\sigma$-bandlimited in QFT sense by variable substitution.  Moreover, all the  conditions of Theorem \ref{main-results-GQPSWFs} are satisfied.
Then we can find a countably infinite set of $\Q$-valued signals ${\{\phi_n(\bx)\}}_{n=1}^{\infty}$ called quaternionic prolate spheroidal wave functions (QPSWFs) and a set of $\Q$-valued numbers
$|\lambda_1|\geq|\lambda_2|\geq\cdot\cdot\cdot$
such that for every $(x_1,x_2)\in \mathbb{R}^2$,
\begin{equation} \label{aaa}
 \lambda_n\phi_n(x_1,x_2)=\frac{1}{2\tau}\int_{[-\tau,\tau]^2}\phi_n(\omega_1,\omega_2)\e^{\qj \sigma x_2\omega_2 \tau}\e^{\qi \sigma x_1\omega_1/ \tau}d\omega_1d\omega_2
\end{equation}
and
\begin{equation} \label{bbbb}
|\lambda_n|^2\phi_n (x_1,x_2)=\int_{{[-\tau,\tau]^2}}\phi_n (y_1,y_2) \frac{\sin \sigma (x_1-y_1)\sin \sigma (x_2-y_2)}{ \sigma(x_1-y_1) \sigma(x_2-y_2)} \,dy_1dy_2 .
\end{equation}
Since $\{E(\w,\bx_{n_1n_2})=\frac{1}{2\tau} \e^{-\qi \omega_1\frac{\pi}{\tau}n_1} \e^{-\qj \omega_2 \frac{\pi}{\tau} n_2}\}_{(n_1,n_2)\in \mathbb{Z}^2}$ is an orthonormal basis of $L^2 (\D,\Q)$, then for any $f\in  \mathcal{H}$,
\begin{equation*}
  \begin{split}
    f(x_1,x_2) =&  \sum_{n_1, n_2} f(\frac{\pi}{\sigma}n_1,\frac{\pi}{\sigma} n_2)\frac{\sin(\sigma x_1-n_1 \pi)\sin(\sigma x_2-n_2 \pi)}{(\sigma x_1-n_1 \pi)(\sigma x_2-n_2 \pi)} \\
     = & \sum_{n=1}^\infty \left(\sum_{m_1, m_2} f(\frac{\pi}{\sigma}m_1,\frac{\pi}{\sigma}m_2)\overline{\phi_n(\frac{\pi}{\sigma}m_1,\frac{\pi}{\sigma}m_2)}\right) \phi_n(x_1,x_2).
  \end{split}
\end{equation*}
From Plancherel theorem, we have
$f\in L^2(\mathbb{R}^2,\Q)$ and
%\begin{equation*}
 % \int_{\mathbb{R}^2}|f(p_1,p_2)|^2dp_1dp_2=\int_D|F(t_1,t_2)|^2dt_1dt_2
%\end{equation*}
\begin{equation*}
 \int_{\mathbb{R}^2} f(x_1,x_2)\overline{g(x_1,x_2)}dx_1dx_2=\frac{\pi^2}{\sigma^2}\int_{[-\tau,\tau]^2} F(\omega_1,\omega_2)\overline{G(\omega_1,\omega_2)}d\omega_1d\omega_2.
\end{equation*}
Notice that the inner product of $\mathcal{H}$ is defined by
\begin{equation*}
  (f,g)_{\mathcal{H}}=(F,G)_H.
\end{equation*}
Thus
\begin{equation*}
  (f,g)_{\mathcal{H}}=\frac{\sigma^2}{\pi^2} \int_{\mathbb{R}^2} f(x_1,x_2)\overline{g(x_1,x_2)}dx_1dx_2
\end{equation*}
and $(\mathcal{H},(\cdot,\cdot)_{\mathcal{H}})$ is a subspace of  $L^2(\mathbb{R}^2,\Q)$.
 By Theorem \ref{reproducing-kernel-Hilbert-spaces}, we have $f(x_1,x_2)$ equals to
 \begin{equation*}
  \int_{\mathbb{R}^2} f(y_1,y_2)\frac{\sin \sigma (x_1-y_1)\sin \sigma (x_2-y_2)}{ \sigma(x_1-y_1) \sigma(x_2-y_2)}dy_1dy_2
\end{equation*}
for any $f(x_1,x_2)\in \mathcal{H}$. Furthermore,
\begin{equation*}
  \displaystyle\sup\left\{\frac{\displaystyle\int_{[-\tau,\tau]^2} |f(x_1,x_2)|^2dx_1dx_2}{\displaystyle\int_{\mathbb{R}^2} |f(x_1,x_2)|^2dx_1dx_2}\right\}=\frac{\sigma^2 \lambda_1^2}{\pi^2}.
\end{equation*}
The extrema is reached if $f(x_1,x_2)=\mathcal{T}\Phi_1(x_1,x_2)=\phi_1(x_1,x_2)$.
\end{ex}

%\begin{rem}
%Suppose that $\{\varphi_{n}(x)\}$ and  $\{\alpha_n\}$ are classical PSWFs and corresponding eigenvalues with parameters $\tau,\sigma$  respectively. Let $\phi_{mn} (x_1,x_2)=\varphi_{m}(x_1)\varphi_{n}(x_2)$ and $\lambda_{mn}=\frac{\sigma^2 }{\pi^2}\alpha_m\alpha_n$. Then $\{(\phi_{mn} (x_1,x_2),\lambda_{mn})\}$  satisfy (\ref{aaa}) and (\ref{bbbb}). Moreover, $\phi_{mn} (x_1,x_2)$ are orthogonal and complete in $L^2(\D,\Q)$ and $\lambda_{00}$ is the maximum of eigenvalues of (\ref{aaa}).
%\end{rem}

% ------------------------------------------------------------------------

\subsection*{Acknowledgements}
 Kit Ian Kou acknowledges financial support from the National Natural Science Foundation of China under Grant (No. 11401606), University of Macau (No. MYRG2015-00058-L2-FST and No. MYRG099(Y1-L2)-FST13-KKI) and the Macao Science and Technology
Development Fund (No. FDCT/094/2011/A and No. FDCT/099/2012/A3).

% ------------------------------------------------------------------------

\begin{thebibliography}{1}


\bibitem{bell1968special}
 Bell, W.W.: Special Functions for Scientists and Engineers. Van Nostrand  (1968)

\bibitem{landau1961prolate1}
  Landau, H.J.,   Pollak, H.O.: Prolate spheroidal wave functions, {F}ourier
  analysis and uncertainty{-}{I}. Bell Syst  Tech J.   40(1), 43--64 (1961)

\bibitem{landau1961prolate}
 Landau, H.J.,   Pollak, H.O.: Prolate spheroidal wave functions, {F}ourier
  analysis and uncertainty{-}{II}. Bell Syst Tech J. 40(1), 65--84  (1961)

\bibitem{landau1962prolate}
 Landau, H.J.,   Pollak, H.O.: Prolate spheroidal wave functions, {F}ourier
  analysis and uncertainty{-}{III}: The dimension of the space of essentially
  time-and band-limited signals. Bell Syst Tech J. 41(4), 1295--1336  (1962)

\bibitem{slepian1964prolate}
Slepian, D.: Prolate spheroidal wave functions, {F}ourier analysis and
  uncertainty{-}{IV}: extensions to many dimensions; generalized prolate
  spheroidal functions. Bell Syst  Tech  J.  43(6), 3009--3057 (1964)

\bibitem{zayed2007generalization}
Zayed, A.I.: A generalization of the prolate spheroidal wave functions. Proc
  Amer  Math  Soc.  135~(7), 2193--2203 (2007)

\bibitem{moumni2014generalization}
Moumni, T.,   Zayed, A.I.: A generalization of the prolate spheroidal wave
  functions with applications to sampling. Integral Transforms  Spec  Funct. 25(6), 433--447  (2014)

\bibitem{walter2003sampling}
 Walter, G.G.,  Shen, X.: Sampling with prolate spheroidal wave functions. J Sampl Theory:
 Signal Image Process.  2(1),  25--52 (2003)

\bibitem{sangwine1996fourier}
 Sangwine, S.J.: Fourier transforms of colour images using quaternion or
  hypercomplex, numbers. Electron  Lett. 32(21), 1979--1980 (1996)

%\bibitem{bihan2003quaternion}
% Bihan NL,  Sangwine SJ. Quaternion principal component analysis of color
%  images, in: Image Processing, 2003. ICIP 2003. Proceedings. 2003
%  International Conference on, Vol.~1, IEEE, 2003, pp. I--809.

\bibitem{ell2007hypercomplex}
 Ell, T.A.,   Sangwine, S.J.: Hypercomplex {F}ourier transforms of color
  images.  IEEE  Trans  Image Process.  16(1), 22--35 (2007)

\bibitem{ell2013quaternion}
 Ell, T.A.: Quaternion {F}ourier transform: Re-tooling image and signal
  processing analysis, in: Quaternion and Clifford Fourier Transforms and
  Wavelets, Springer,  pp. 3--14 (2013)

\bibitem{brackx1982clifford}
 Brackx, F.,  Delanghe, R.,  Sommen, F.: Clifford Analysis. Vol. 76, Pitman Books
  Limited  (1982)

\bibitem{ghiloni2013continuous}
 Ghiloni, R.,  Moretti, V.,  Perotti, A.: Continuous slice functional calculus in
  quaternionic  Hilbert spaces. Rev  Math  Phys.  25(04), 1350006 (2013)

\bibitem{fashandi2013compact}
 Fashandi, F.: Compact operators on quaternionic Hilbert spaces. Facta Univ
  Ser  Math  Inform.  28(3), 249--256 (2013)

\bibitem{farenick2003spectral}
 Farenick, D.R.,  Pidkowich, B.A.: The spectral theorem in quaternions. Linear
  Algebra Appl. 371, 75--102 (2003)

\bibitem{ghiloni2014spectral}
 Ghiloni, R.,  Moretti, V.,  Perotti, A.: Spectral properties of compact normal
  quaternionic operators. in: Hypercomplex Analysis: New Perspectives and
  Applications, Springer, pp. 133--143 ( 2014)

\bibitem{Alpay2016unbounded}
Alpay, D., Colombo, F., Kimsey, D.P.: The spectral theorem for quaternionic unbounded normal operators based on the $S$-spectrum. J. Math. Phys. 57, 023503 (2016)

\bibitem{Alpay2016unitary}
Alpay, D., Colombo, F., Kimsey, D.P., Sabadini, I.: The Spectral theorem for unitary operators based on the $S$-spectrum. Milan J. Math.   84(1), 41--61 (2016)


\bibitem{Colombo2016Schatten}
  Colombo, F., Gantner, J., Janssens, T.:  Schatten class and Berezin transform of quaternionic linear operators. Math. Meth. Appl. Sci. 39(2016), 5582--5606 (2016)



\bibitem{rudin1987real}
 Rudin, W.: Real and Complex Analysis. Tata McGraw-Hill Education (1987)

\bibitem{bollabas1999linear}
 Bollab{\'a}s, B.: Linear Analysis, An Introductory Course. Cambridge University
  Press  (1999)

\bibitem{Colombo2011Noncommutative}
  Colombo, F., Sabadini, I. Struppa, D.C.: Noncommutative Functional Calculus: Theory and Applications of Slice Hyperholomorphic Functions. Birkh{\"a}user, Basel (2011)

%\bibitem{tobar2014quaternion}
% Tobar, F.,   Mandic,  D.P.: Quaternion reproducing kernel {H}ilbert spaces:
%  Existence and uniqueness conditions.  IEEE  Trans  Inf  Theory.  60(9), 5736--5749 (2014)
\bibitem{Alpay2016SliceSchur}
Alpay, D., Colombo, F., Sabadini, I.: Slice Hyperholomorphic Schur Analysis, Operator Theory: Advances and Application. Birkh{\"a}user (2016)


\bibitem{saitoh1983hilbert}
 Saitoh, S.: Hilbert spaces induced by {H}ilbert space valued functions. Proc
  Amer  Math  Soc.  89(1), 74--78 (1983)

\end{thebibliography}
\end{document}